\documentclass{article}

\usepackage{t1enc}
\usepackage{times}
\usepackage{amsmath}
\usepackage{amsfonts}
\usepackage{pstricks}
\usepackage{graphicx,float}
\usepackage{longtable}
\newtheorem{theorem}{Theorem}
\newtheorem{definition}{Definition}

\title{A preconditioning strategy for microwave susceptibility in ferromagnets}

\author{St\'ephane Labb\'e\thanks{Universit\'e Paris Sud, Laboratoire de Math\'ematique, B\^at. 425, 91405 Orsay Cedex,  33-(0)1-69-15-60-42,  33-(0)1-69-15-67-18, stephane.labbe@math.u-psud.fr}}

\begin{document}

\maketitle

\begin{abstract}
3D numerical simulations of ferromagnetic materials can be compared with experimental results via 
microwave susceptibility. In this paper, an optimised computation of this microwave susceptibility 
for large meshes is proposed. The microwave susceptibility is obtained by linearisation of the Landau and Lifchitz
equations near equilibrium states and the linear systems to be solved are very ill-conditionned.
Solutions are computed using the Conjugate Gradient method for the Normal equation (CGN Method).
An efficient preconditioner is developed consisting of a projection and  an approximation of an ``exact'' preconditioner  
in the set of circulant matrices. Control of the condition number due to the preconditioning 
and  evolution of the singular value decomposition are shown in the results.
\end{abstract}

\section{Introduction}
Ferromagnetic simulation via the micromagnetic model is a real-life computational challenge. Ferromagnetic materials are used in numerous
applications such as radar protection, magnetic recording or micro electronics. In these applications, the magnetic objects studied are micro
or nano-objects which are difficult and expensive to craft. Thus, one of the optimisation solutions, for the shape and composition of such
particles, is numeric simulation. The first step in this type of simulation is to compute the dynamic of the magnetisation and the equilibrium
states. However, a direct comparison of the results with experiments is impossible for 3D particles. The main comparison tool is
microwave susceptibility as the resonnance numerical curves can be compared with the physical experiments. At that point several difficulties are
encountered. The main one is managing a large number of degrees of freedom. This is  required to compute interesting configurations with sufficient accuracy.

In this article, we use the micromagnetism model in order to model the magnetisation behaviour in ferromagnetic materials. This model
is a mesoscopic model, ie. a model valid for a scale between the one used for microscopic Maxwell equations and the scale of 
classic macroscopic Maxwell equations. In this model, magnetisation does not linearly depend on magnetic excitation but is controlled
by a non-linear system: the Landau-Lischitz equation (\ref{LL}). This model was introduced by Brown \cite{Brown:Microm,Brown:magnet}.

There are two ways to obtain the equilibrium states. The first by energy minimisation (\cite{Alouges:habilitation,Bagneres:DeuxD,Schabes.Bertram:Magne}), the second 
by relaxation of the dynamic system (\cite{Miltat:Dynamical,Labbe.Bertin:Microwave}). The main advantage of the dynamical
approach is to compute an equilibrium state linked to  given initial data by a life-like dynamic process; then, we can apply
dynamical treatments, via the external field, in order to find specific equilibrium states. 

Computation of the microwave susceptibility can be performed by two main methods: the Harmonic Direct Computation and the Fourier Transform
Method. The first method is based upon the use of a linearised version of the evolution equation pertubated by a time harmonic external
field. The second is based upon the injection of an harmonic perturbation. The Fourier method implies the resolution of a time dependant
problem that is quite ill-conditionned for low frequencies (the time step ensuring that the convergence vanishes swiftly when the frequency
decreases) but the linearisation methods permit the range of frequencies used in the applications to be attained.

\section{The microwave susceptibility problem}
\subsection{The linearisation}
In this problem, we are interested in computing the microwave response of a ferromagnetic system to an external
harmonic exitation. We consider that the ferromagnetic material is homogeneous and contained in a $\mathcal{C}^1$-class
piecewize domain of $\mathbb{R}^3$ denoted $\Omega$. Then, we study the evolution of the magnetisation field in 
the neighbourhood of an equilibrium state of the dynamic equation. This equation, in the micromagnetism model 
\cite{Brown:Microm}, is given by the Landau-Lifchitz system: find $m$ in $\widetilde{H}^1([0,T]\times\mathbb{R}^3,\Omega;\mathbb{R}^3)$ 
$=\{m\in L^2(\mathbb{R}^3;\mathbb{R}^3)\vert \forall t\in[0,T],\;m_{\vert\Omega}\in H^1(\Omega;\mathbb{R}^3) \mbox{ and } m\equiv 0 \mbox{ in } \mathbb{R}\backslash\Omega\}$ such that
\begin{equation}
\left\{
\begin{array}{l}
\frac{\partial m}{\partial t}=f(m,h_{ext})=-m\wedge (H(m)+\ell)-\alpha\ m\wedge(m\wedge (H(m)+\ell)),\in (0,T]\times\Omega,\\
m(x,0)=m_0(x),\;\forall t\in\Omega.
\end{array}
\right.
\label{LL}
\end{equation}
where $H$ is a linear operator, from $\widetilde{H}^1([0,T]\times\mathbb{R}^3,\Omega;\mathbb{R}^3)$ into $H^{-1}(\mathbb{R}^3;\mathbb{R}^3)$, 
$\ell$ the external magnetic field (independent of the magnetisation and element of $L^{\infty}([0,T]\times\mathbb{R}^3;\mathbb{R}^3)$,
$\alpha$ the damping factor (a strictly positive real) and $m_0$ is a given element of $\widetilde{S}^2(\Omega)$=$\{m\in\widetilde{H}^1(\mathbb{R}^3,\Omega;\mathbb{R}^3)\; \vert\; \vert m_{\vert\Omega}\vert=1,\mbox{ a.e. in } \Omega\}$. In this model, we can see that the local module of the magnetisation is naturally preserved.
In this article, we define $H$ as follows: $\forall m\in \widetilde{H}^1([0,T]\times\mathbb{R}^3,\Omega;\mathbb{R}^3)$
$$
H(m)=A\triangle m + H_d(m) + K (m - (m.u)u)  
$$
where $A$ and $K$ positive real constants and $u$ is an element of $\widetilde{H}^1([0,T]\times\mathbb{R}^3,\Omega;S^2)$ ($S^2$ designates
the unit sphere). The operator $H_d$ is defined in the sense of distributions on $\mathbb{R}^3$ by 
$$
\left\{
\begin{array}{l}
\mbox{rot}(H_d(m))=0,\\
\mbox{div}(H_d(m))=-\mbox{div}(m).
\end{array}
\right.
$$

Now, let us define the equilibrium states of the system (\ref{LL})
\begin{definition}
For a given $\ell$ in $L^{\infty}(\mathbb{R}^3;\mathbb{R}^3)$ (independent of time), a magnetisation state $m_\ell$, in $\widetilde{H}^1(\mathbb{R}^3,\Omega;\mathbb{R}^3)$
is an equilibrium state if, and only if,
$$
f(m_\ell,\ell)=0,\;\mbox{a.e. in }\Omega.
$$
\end{definition}
Then, for a given equilibrium state $m_\ell$, associated to an external state $\ell$, we define the microwave 
susceptibility
\begin{definition}
For a given equilibrium state $m_\ell$, associated to an external field $\ell$, we denote a susceptibility tensor
of the order 3 complex matrices $\chi(\ell)$ defined by
$$
(\chi(\ell))_{l,k}=-\frac{1}{2\ T}(\lambda_k,m_l)_{0,\Omega},\;\forall (l,k)\in\{1,2,3\}^2,
$$ 
with $\lambda_k=\zeta_k e^{i\omega t}$ and $\zeta_k$ is a contant vector of $\mathbb{R}^3$. Furthermore, we suppose
that $(\zeta_k)_{k\in\{1,2,3\}}$ is an orthogonal basis of $\mathbb{R}^3$. Then, for all $k$ in $\{1,2,3\}$, $m_k$
is a solution of (\ref{LL}) for the external field $\lambda_k+\ell$ and the intial data $m_0=m_\ell$.
\end{definition}
Formally, if the excitation $\zeta_k$ is sufficiently small, then the magnetisation responses will be also
small and we can define this response for every $k$ in $\{1,2,3\}$ by
$$
m_k-m_\ell=\mu_k\ e^{i\omega t},
$$
with $\mu_k \in \widetilde{H}^1(\mathbb{R}^3,\Omega;\mathbb{C}^3)$. In the following we suppose that
$\zeta_k$ and $\mu_k$ are of the same order.

Then, if we re-write the system (\ref{LL}) verified by $m_k$, the linearised equation gives
\begin{equation}
(i\omega - D_{1,\ell}\circ h - D_{2,\ell})(\mu_k)=D_{1,\ell}(\zeta_k)
\label{LLL}
\end{equation}
where, for all $w$ in $L^\infty(\mathbb{R}^3;\mathbb{R}^3)$, we set
$$
D_{1,\ell}(w)=-m_\ell\wedge w-\alpha\ m_\ell\wedge(m_\ell\wedge w),
$$
$$
D_{2,\ell}(w)=(H(m_\ell)+\ell)\wedge w + \alpha\ m_\ell \wedge(w\wedge (H(m_\ell) + \ell))
$$
\subsection{The discretisation of the linearised equation}
In order to discretise the equation, we consider a monolith $K(\Omega)$ such that $\Omega\subset K(\Omega)$.
Ideally, this monolith is the smaller containing $\Omega$. Then, $K(\Omega)$ is discretised using a regular
cubic mesh of cells $(\Omega_i)_{i\in N_h}$ where $h$ is the length of a cell and $N_h$ is the set of the 
indices. We set $\Omega_h=\bigcup_{i\in N_{int,h}}\Omega_i$ where $N_{int,h}\subset N_h$ is the set of indices such 
that, for every $i$ in $N_{int,h}$, $\Omega_i\cap\Omega\ne \emptyset$.

Then, we choose as a discrete space for all euclidian space $F$:
$$
W_h(F)=\{u\in L^2(\mathbb{R}^3;F)\vert u\equiv 0 \mbox{ in }\mathbb{R}^3\backslash K(\Omega)\mbox{ and }\forall i\in N_h,\;u_{\vert\Omega_i} \mbox{ is a constant}\},
$$
for each $u$ in $W_h$, we set: $\forall i\in N_h,\;u_i=u_{\vert \Omega_i}$.
We choose the $L^2$ scalar product on $\mathbb{R}^3$ as the scalar product on $W_h$, we denote it  $(u,v)_{0,\Omega}$ for
all $u$,$v$ in $L^2(\mathbb{R}^3;F)$.
Then, setting 
$$
\begin{array}{lcl}
                     &   P_h        &\\
L^2(\mathbb{R}^3;F)  & \longrightarrow  & W_h(F)\\
\displaystyle u                    & \longmapsto & \displaystyle  P_h(u)=\sum_{i\in N_h} \left(\frac{\mbox{\bf 1}_{i}}{h^3}\int_{\Omega_i}u\ dx\right)\\
\end{array}
$$
where $\mbox{\bf 1}_{i}$ is defined for $x$ in $\mathbb{R}^3$ by $\mbox{\bf 1}_{i}(x)=1$ if $x$ belongs to $\Omega_i$, $\mbox{\bf 1}_{i}(x)=0$ otherwise. $P^\star_h$
designates the canonical injection of $W_h(F)$ onto $L^2(\mathbb{R}^3;F)$.

These definitions lead to the following formulas for the discrete magnetic contributions:
$$
H_{a,h} = P_h\circ H_{a}\circ P^\star_h,
$$ 
and 
$$
H_{d,h} = P_h\circ H_{d}\circ P^\star_h,
$$ 
the analysis of $H_{a,h}$ is straightforward. On the other hand, the analysis of $H_{d,h}$ is not direct, in particular,
it has been demonstrated that this discretisation preserves the main properties of the demagnetisation operator $H_d$ ($H_d$
is a projection operator), and a lower estimate of its lower eigenvalue is given. Furthermore, the computation of this operator
is very expensive: the discrete matrix is a full matrix. Then, to optimise its computation, we choose to use a regular
cubic mesh which ensure a specific structure for the discrete operator. This block-Toeplitz structure enables us to reduce the
storage of the matrix from $\#(N_h)^2$ to $O(\#(N_h))$ and the computation cost from $\#(N_h)^2$ to $O(\#(N_h)\log(\#(N_h)))$.
For complete analysis of the discretisation of $H_d$, see \cite{labbe:fast}.
The Laplacian operator is discretised using the classical $7$ point scheme, the discretised operator is designated in the
following by $\triangle_h$. The total discretised magnetic field is then defined by
$$
H_h(m)=A\triangle_h m+H_{d,h}(m)+H_{a,h}(m).
$$

Then, for a given external field $\ell$ in $W_h(\mathbb{R}^3)$, we set $m_{h,\ell}$, element of $W_h(S^2)$, the equilibrium state of 
the discretised version of (\ref{LL}). This state is obtained using an explicit time discretisation combined with an optimisation
of time which ensures its stability (see \cite{labbe:mumag,Halpern.Labbe:Modelisation,Labbe.Bertin:Microwave}). 
This equilibrium state, as seen previously, is such that: $\forall i\in N_h,\;\exists \beta_i\le 0$ and
$$
H_h(m_{h,\ell}) \mbox{\bf 1}_i=\beta_i m_{h,\ell} \mbox{\bf 1}_i,
$$
we set $H_h(m_{h,\ell})=B_\ell(m_{h,\ell})$ where $B_\ell$ is a diagonal operator.
Knowing an equilibrium state for the discretised sytem, we can define the linearised discrete system: $\forall \omega \in \mathbb{R}^+_*$,

\begin{equation}
(i\omega - D_{1,h,\ell}(H_h-B_\ell))\mu_h=D_{1,h,\ell}\zeta_h
\end{equation}
where $D_{1,h,\ell}$ is the operator $D_{1,\ell}$ built for the $m_{h,\ell}$ equilibrium state.

Then, for each element $u$ of $W_h(\mathbb{R}^3)$, we associate a unique element $U$ of $\mathbb{R}^{3 \#(N_h)}$ defined
by
$$
\forall i\in N_h,\;U_i\in\mathbb{R}^3\mbox{ and }U_i=\frac{1}{h^3}\int_{\Omega_i}u(x)\ dx.
$$ 
Using this bijection between  $W_h(\mathbb{R}^3)$ and $\mathbb{R}^{3 \#(N_h)}$, we can write a matricial version of the 
linearised discrete version of (\ref{LL}): find $U_k$ in $\mathbb{R}^{3 \#(N_h)}$ such that, for a given $Y_k$ built
on $\zeta_k$ we have
\begin{equation}
M_\omega U_k = D\ Y_k,\label{LLLh}
\end{equation}
where, for every $U$ in $\mathbb{R}^{3 \#(N_h)}$, for every $i$ in $N_h$
$$
(M_\omega U)_i = \frac{1}{h^3}\int_{\Omega_i}\left((i\omega - D_{1,h,\ell}(H_h-B_\ell))\left(\sum_{i\in N_h} U_i\mbox{\bf 1}_i\right)\right)\ dx,
$$
and 
$$
(D\ Y_k)_i = \frac{1}{h^3}\int_{\Omega_i}\left(D_{1,h,\ell}\left(\sum_{i\in N_h} Y_{k,i}\mbox{\bf 1}_i\right)\right)\ dx,
$$
For use in the remainder of this paper for every $U$ in $\mathbb{R}^{3 \#(N_h)}$ we set:
$$
D H U = - M_\omega U + i\omega U
$$
then $H$ is the matrix associated to the discrete operator $H_h - B_\ell$.

\subsection{Some properties of the discrete system (\ref{LLLh})}
We set $M_\ell$, the element of $\mathbb{R}^{3 \#(N_h)}$ associated to $m_{h,\ell}$. Let us define $[m_{h,\ell}]^\perp$ by
$$
[m_{h,\ell}]^\perp = \{W\in \mathbb{C}^{3 \#(N_h)} \vert \forall i\in N_h,\;M_{\ell,i}.W_i=0\},
$$
and we designate by $P_\ell^\perp$ the projection from $\mathbb{C}^{3 \#(N_h)}$ into $m_{h,\ell}$. Then we can
demonstrate:
\begin{theorem}
For every $Y$ in $\mathbb{R}^{3 \#(N_h)}$ and for every $\omega$ strictly positive, the system (\ref{LLLh}) is regular and its solution 
is in an element of $[m_{h,\ell}]^\perp$.\label{th1}
\end{theorem}
{\bf Proof:}
If $U$ is the solution of (\ref{LLLh}), then we have
$$
i\omega U = D (Y_k + H U), 
$$
knowing that $D$ sends elements of $\mathbb{C}^{3 \#(N_h)}$ in $[m_{h,\ell}]^\perp$, we conclude that $U$ is also an element
of $[m_{h,\ell}]^\perp$.

Then, considering $V$ in $[m_{h,\ell}]^\perp$, due to the structure of $H$, we have $H V$ as an element of $[m_{h,\ell}]^\perp$.
Each diagonal block (3$\times$3) has 0,$\alpha+i$ and $\alpha-i$ as eigenvalues. Knowing that the eigenvalues of $H$ (symetric matrix)
are real, we deduce that the eigenvalues of $D H$ are complex numbers of non vanishing real parts unless the eigenvalueis null. Then,
the eigenvalues of $M_\omega$ can not vanish.
\begin{flushright}$\Box$\end{flushright}

The conditioning number of the matrix $M_\omega$ can be estimated
\begin{theorem}
For every $\omega$ real strictly positive, we have
$$
\mbox{cond}(M_\omega)\le \sqrt{\frac{\omega^2+(1+\alpha^2)(1+\frac{1}{h^3})(\frac{A}{h^2}+1+K)}{\omega^2}}
$$\label{th2}
\end{theorem}
{\bf Proof:}
This theorem is proved using the Courant-Fisher theorem for hermitian matrices which provides formulae for the highest and
lowest eigenvalues. The proof is then classical and uses the fact that $D^* D$ is the projection matrix on $[m_{h,\ell}]^\perp$
multiplied by $(1+\alpha^2)$.
\begin{flushright}$\Box$\end{flushright}
We notice that the conditioning number $\mbox{cond}(M_\omega)$ bahaves as expected when $\omega$ tends to infinity:
$$
\lim_{\omega\rightarrow \infty} \mbox{cond}(M_\omega) = 1.
$$
Here, the fact that $\omega$ grows to infinity means that it dominates $\frac{1}{h^2}$. Now, if we consider that $\omega$
is fixed, the behaviour of $\mbox{cond}(M_\omega)$ shows that the system is ill-conditioned
$$
\lim_{h\rightarrow 0} \mbox{cond}(M_\omega) = \infty.
$$
Thus, the pre-conditioning of the system is essential. In fact, the most interesting part of the spectrum of susceptibility
for numerous applications is the low frequency part.

\section{The precontioning strategy}
\subsection{Choice of the inversion method}
In order to solve system (\ref{LLLh}), we chose an iterative method; this choice is conditioned by the fact that the matrices 
considered are non-symmetric full matrices and the order of the systems to solve is great (up to $10^6$). Three main iterative
methods are used commonly to solve non symmetric systems:
\begin{itemize}
\item{the normal conjugate gradient (CNG),}
\item{the generalised minimal residual method (GMRES),}
\item{the conjugate gradient squared (CGS).}
\end{itemize}
As shown in the article of Nachtigal, Reddy and Trefethen \cite{Nachtigal:How.fast}, none of this three methods could be considered
as a cure-all for all non-symmetric systems. As the convergence quality of CGS and GMRES is influenced by eignevalue clustering
of the system matrix, CNG method convergence depends on singular value clustering. As the preconditionning strategy 
presented in this article is based upon the amelioration of the singular value clustering, we chose, of course, the CNG method. Furthermore,
tests not presented in this article show that the CNG method seems to be more adaptated for this type of system, even if not
preconditioned.
\subsection{An example of singular value repartition and of CNG convergence rate}
In the remainder of this paper, we have chosen to illustrate the results presented using a plain example. This example has been chosen for
the low order, 192, of its system which facilitates the visualisation (done with Matlab). The mesh chosen is a $4\times 4\times 4$ regular
cubic mesh of a cubic domain. We set it in the dimensionless system $A=0,88\ 10^{-10}$, $K=0,57\ 10^{-2}$, $\alpha=0,5$ and the cube length is equal 
to $10^{-6}$. For this bench, we would want to choose $\omega$ between $\omega_{min}=0,452\ 10^3$ Hz and $\omega_{max}=0,452\ 10^{5}$ Hz.
\begin{figure}[h]
\begin{center}
\includegraphics[width=7cm]{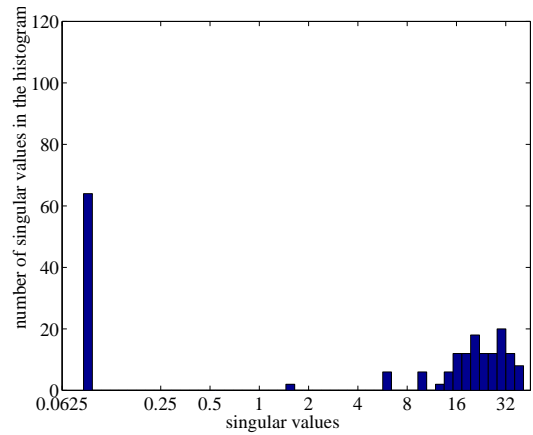}
\end{center}
\caption{Singular value decompostion for the non-preconditioned system.}
\label{fig1}
\end{figure}
In Fig. \ref{fig2} the error evolution for the CNG is shown. Here we have chosen a final error criteria of $10^{-5}$. With no preconditioning,
the system converges in 56 iterations for $\omega_{min}$ and 48 iterations for $\omega_{max}$, the precontioning number is almost equal
to 7500 (slight variations between  $\omega_{min}$ and $\omega_{max}$).
\begin{figure}[h]
\begin{center}
\includegraphics[width=7cm]{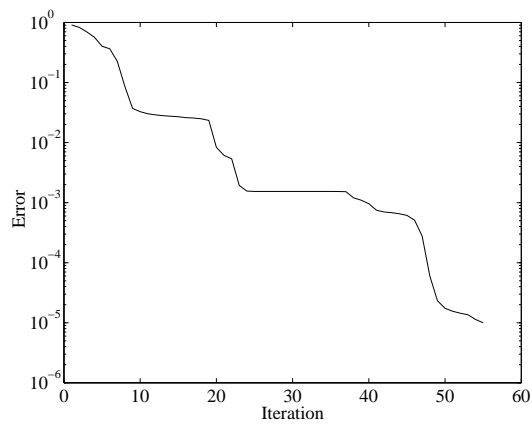}
\end{center}
\caption{Residue of the CNG for $\omega_{min}$.}
\label{fig2}
\end{figure}
\subsection{The preconditioning strategy}
We have three main goals to build the preconditioner:
\begin{itemize}
\item{to use the known properties of the system,}
\item{decrease the conditioning number sensitivity to the mesh size,}
\item{build a cheap preconditioner (memory size and computational cost).}
\end{itemize}
The first point is taken into account by using the result presented in Theorem~\ref{th1}: the right side of
the preconditioner will be a projection on $[m_{h,\ell}]^\perp$. This first step of projection eliminates the
cluster of singular values near $0$ and ensures a convergence in 48 iterations for $\omega_{min}$ and of 27
iterations for $\omega_{max}$.

%\begin{figure}[h]
%\begin{center}
%\includegraphics[width=7cm]{svdNP}
%\end{center}
%\caption{Singular value decompostion for the right preconditionned system.}
%\label{fig1}
%\end{figure}

\subsection{The ``exact'' preconditioner}
As a first stage, we would want to build a symmetric left precontioner. The non symmetry of $M_\omega$
comes from the operator $D_{1,h,\ell}$. In fact, we have
$$
D_{1,h,\ell}\ w=-m_\ell\wedge w+\alpha P^\perp_\ell w,
$$
the first part of the operator is a rotation, and the second part a projection. Then, it is possible
to prove that $D_{1,h,\ell}$ does not have a main influence on the singular value decomposition. This means that we
may choose a left preconditioner $M_{g,\omega}$ built on the operator
$$
i\omega - \alpha (H_h - B_\ell),
$$
That is to say, if we set $H$ the matrix built on the operator $H_h - B_\ell$
$$
M_{g,\omega}=i\omega Id - \alpha H.
$$
In the sequel, even if we do not write the projection to lighten the notations, we consider that the system is right
preconditioned by $P^\perp_\ell$. Then, we prove the following theorem
\begin{theorem}
For each $\omega$ and $h$ strictly positive, we have 
$$
\mbox{cond}(M^{-1}_{g,\omega}M_\omega)\le \sqrt{1+g(h,\omega)\left[2+g(h,\omega)^2\right]}
$$
where
$$
g(h,\omega)=\sqrt{\frac{(2+\alpha^2)(1+\frac{1}{h^3})^2(\frac{A}{h^2}+1+K)^2}{\omega^2+(1+\frac{1}{h^3})^2(\frac{A}{h^2}+1+K)^2}}
$$
\label{th3}
\end{theorem}
{\bf Proof:}
In the space $[m_{h,\ell}]^\perp$, we have: $D=R+\alpha Id$, where $Id$ is the eye matrix on space $[m_{h,\ell}]^\perp$
and $R$ is the matrix associated to the operator $-m_\ell \wedge$. Moreover, for every $U$ in $[m_{h,\ell}]^\perp$, we have
$$
(H U,m_{h,\ell})=(U, H m_{h,\ell})=0,
$$
this implies, by breaking off of the elements of $[m_{h,\ell}]^\perp$, that $H U$ is an element of $[m_{h,\ell}]^\perp$.
We remark also that by working in  $[m_{h,\ell}]^\perp$, we have
$$
R^2=-Id.
$$

Then, we have
$$
\begin{array}{lcl}
M^{-1}_{g,\omega}M_\omega&=&(i\omega H^{-1} - Id) H^{-1} N H,\\
\end{array}
$$
where 
$$
N_h=-R\ H + (1 - \alpha) H=N H.
$$
Then, for every $V$ in $[m_{h,\ell}]^\perp$, we have the following relation
$$
\begin{array}{lcl}
M^{-1}_{g,\omega}M_\omega V.M^{-1}_{g,\omega}M_\omega V &=& \Vert V\Vert^2 + 2\mathcal{R}[V.M^{-1}_{g,\omega}\ N_h V] + \Vert M^{-1}_{g,\omega}N_h\ V\Vert^2.
\end{array}
$$
We designate as $\mathcal{R}[z]$ the real part of a complex number $z$.

Furthermore, we have the following estimations: 

$$
\Vert (H^{-1} R\ H)^2\Vert = \Vert (H^{-1} R^2\ H)^2\Vert = 1,
$$
this implies that $\Vert H^{-1} R\ H\Vert=1$. So, using the fact than
$$
H^{-1} N H=-H^{-1} R H +(1-\alpha) Id,
$$
we have
$$
\Vert H^{-1} N\ H\Vert \le \sqrt{2 +\alpha^2}
$$
and

$$
\Vert M^{-1}_{g,\omega}\Vert = \max_{j\in N_h^\perp}(\frac{\omega^2}{\lambda_j^2}+1)^{-1} \le g(h,\omega).
$$
where $\lambda_j$ is the eigenvalues of the matrix $H$ in $[m_{\ell,h}]^\perp$ and $N_h^\perp$ is the set of indeces
of $[m_{\ell,h}]^\perp$.

Then, using  the lowest eigenvalue  controlled by projection part of the precondioner we conclude the proof of the
Theorem.

Finally, we have the good behaviour of the preconditionned system when the mesh length $h$ tends to $0$:
$$
\lim_{h\rightarrow 0}(M^{-1}_{g,\omega}M_\omega)\le 1+\sqrt{2+\alpha^2},
$$
This version of the preconditioner gives excellent control of the conditionning number but needs
the inversion of a full matrix. This leads us to the second stage in which we will replace the complete operator
$H_h$ by its laplacian part. 

\subsection{Preconditioning by the Laplacian component: the direct approach}
\begin{figure}[h]
\begin{center}
\includegraphics[width=7cm]{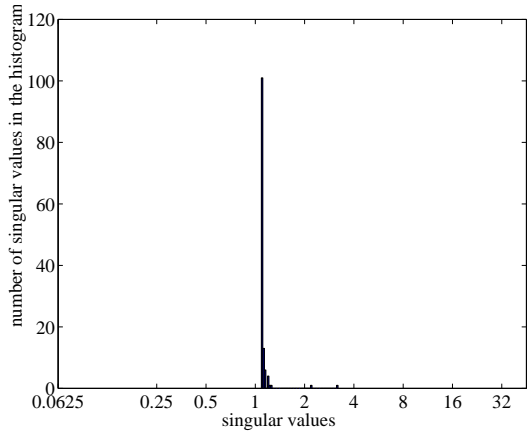}
\end{center}
\caption{Singular value decomposition of the system preconditioned by $M_{g,\omega,\triangle}$.}
\label{fig3}
\end{figure}

\begin{figure}[h]
\begin{center}
\includegraphics[width=7cm]{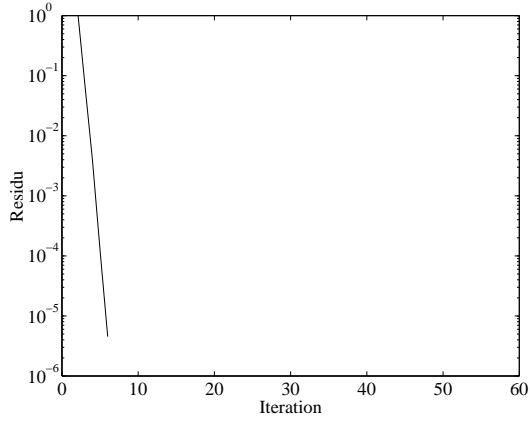}
\end{center}
\caption{Residue of the system preconditionned by $M_{g,\omega,\triangle}$.}
\label{fig4}
\end{figure}

The Laplacian part of $H_h$ is the most punitive part of the matrix $M_\omega$ in terms of preconditioning. The idea 
in this section is to develop an approximate conditioner $M_{g,\omega,\triangle}$ built on the operator
$$
i\omega-A\triangle_h-B.
$$
The matrix $M_{g,\omega,\triangle}$ is a band matrix which could be more easily handled than $M_{g,\omega}$, the
earlier version of the preconditioner. This approximation of the preconditioner $M_{g,\omega}$ will be all the more
accurate as the norms of the operators $H_d$ and $H_a$ are dominated by $\frac{A}{h^2}$. As seen in Fig. \ref{fig3}, the
clustering of the singular value decomposition obtained for the system preconditioned by $M_{g,\omega,\triangle}$ is good.
The convergence of the CNG algorithm using this method is very good: 7 iterations for $\omega_{min}$ and 9
iteration for $\omega_{max}$ (see Fig. \ref{fig4}).

\subsection{The approximated preconditioner}
Nevertheless, the use of the pre-conditioner $M_{g,\omega,\triangle}$ stays expansive. The solution is to build an easily 
invertible approximation of $M_{g,\omega,\triangle}$. Here we will use here the work of \cite{Tyr:optim}. The idea is to project
the matrix $M_{g,\omega,\triangle}$ into the circulant matrix space in the sense of the Froebenuis norm. 

Given a circulant matrix $C$ $n$ by $n$ on $\mathbb{C}$ generated by $c$, vector of $\mathbb{C}^n$, we have
$$
\forall (i,j)\in \{1,...,n\}^2 \mbox{ and } p\in \{1,...,n\},\;C_{i,j}=c_{p}\;\mbox{ if } j-i=p-1 \;\mbox{ or } j-i=n-p1.
$$
Then, as shown in \cite{Tyr:optim}, for every matrix $M$  $n$ by $n$ on $\mathbb{C}$, the projection $C$ of $M$ on 
the space of the circulant matrices of order $n$ is generated by the vector $c$ given by
$$
\forall p\in \{1,...,n\},\;c_p=\frac{1}{n}\left(\sum_{l=1}^{n-p+1} M_{l,l-p+1}+\sum_{l=n-p+2}^{n} M_{l,l+n-p+1}\right) 
$$

The three dimensional projection is more complex but the main idea is contained in the one-dimensional projection.

When the circulant approximation matrix is built, the inversion is performed in the Fourier space (the matrix produced
is block-diagonal 3$\times$3 in Fourier space), then the precondioning is of complexity $O(N\log(N))$ for each iteration
of the inversion method. The other main advantage of the method is that the storage is reduced to $O(N)$.

In this section, we have to keep in mind that the structure is a three dimensional one: the considered matrices are 3 level block
matrices. This implies that the projection must be performed on the 3 levels block circulant matrices.

In the small example presented to illustrate the paper, the system precontioned by the approximated preconditioner 
converges in 30 iterations for the smaller frequency and 26 iterations for the highest (see Fig. \ref{fig6}). The convergence
curve is very good in the sense that the slope is quasi-constant. This point is quite important: susceptibility computations 
do not need high numerical accuracy. Effectively, the results obtained will be compared to experimental results
for which the error is quite important. This comes form the fact that the samples used for experiments are far to be perfect
and that the measurement tools do not have very high precision for this type of experiment.

\begin{figure}[h]
\begin{center}
\includegraphics[width=7cm]{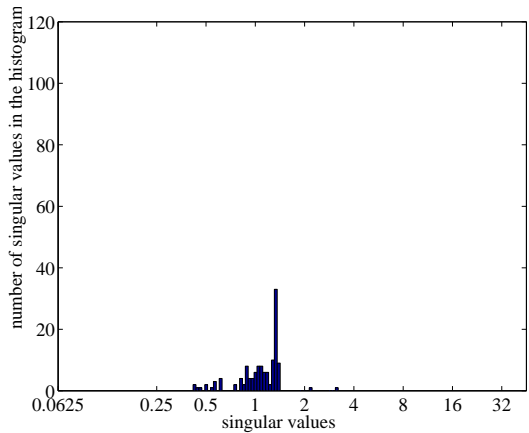}
\end{center}
\caption{Singular value decomposition of the system preconditioned by the circulant approximation of $M_{g,\omega,\triangle}$.}
\label{fig5}
\end{figure}

\begin{figure}[h]
\begin{center}
\includegraphics[width=7cm]{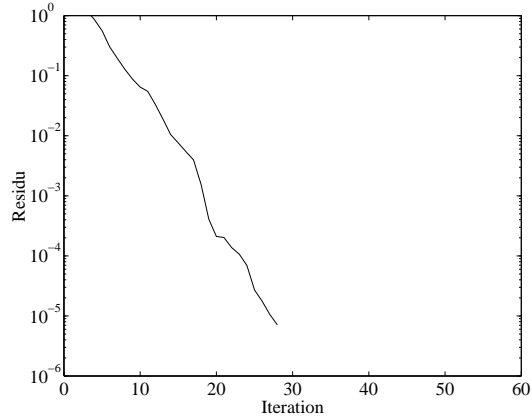}
\end{center}
\caption{Residue of the system preconditioned by the circulant approximation $M_{g,\omega,\triangle}$.}
\label{fig6}
\end{figure}

\section{Numerical simulations}
We present here the number of iterations for the simulation of a ferromagnetic dot. This dot is meshed
by a regular grid, size 64$\times$64$\times$32. In this monolith a cylinder with the axis $z$ and a
circular basis (32 cells for the $z$ direction and 64$\times$64 for the others)is included. The total number of degrees
of freedom is 393216. The results shown here have been computed on the parallel machines of ONERA and
Dassault Aviation.\\

\subsection{Parallel implementation}
There are two possible levels of parallelisation for this problem: local parallelisation for computations of each
iteration and global parallelisation of the frequency computations.

The global parallelisation is a repartition of each frequency computation through the processors. A main process
distributes the computation to each processor such that each processor is always occupied. This part is implemented 
using MPI.

The local implementation, not used for the results presented here, is the parallelisation of the total magnetic field
over the domain. In this computation, one part is more expensive than the others: the demagnetisation field. In fact,
the computation of demagnetisation is accelerated by using its Toeplitz structure (see \cite{labbe:fast}). This 
computation strategy uses 3 dimensionnal FFT intensively. To enhance the performance, we have to parallelise the
FFT computation. To do so, we have chosen to use OPEN-MP. This choice avoids the transposition of the data via the cluster
that must be performed while using a distributed memory system. The results are very satisfying: for a cubic structure
and sufficient number of cells (for instance a $32\times 32\times 32$ mesh), the computation time of FFT is divided
by $1.9$ on a node of two processors. 

\subsection{Description of the benchmark statistics and results}

The aim is to compute the susceptibility of a cylinder of permalloy (see for example \cite{Boust:High} for this
type of results). The parameters of the material are the following:
\begin{center}
\begin{tabular}{|l|c|}
\hline
Parameter & Value \\
\hline
$A$      & 0.17875 $10^{-11}$\\
$\alpha$ & 0.05\\
\hline
\end{tabular}
\end{center}

In the following table \ref{res}, we give the number of iterations for 
directions $x$ and $y$. The direction $z$ in this computation
is omitted because there is no resonnance in this direction.

The computation has been carried out on a node composed of  8 Power4 IBM (1.1GHz) with 16 GO of Ram. An iteration
takes almost 24 seconds, the complete computation took 36 hours.

\begin{longtable}{|l|l|l|l|l|}
\hline
$\omega/1,356$  (Hz) &  iterations  for $x$ & error &  iterations i for $y$ & error \\
\hline\endhead
\hline
\endfoot
3.00 $10^{5}$  &  57  &  4.92 $10^{-2}$  &  127  &  4.98 $10^{-2}$\\
2.73 $10^{5}$  &  100  &  4.85 $10^{-2}$  &  137  &  4.91 $10^{-2}$\\
2.49 $10^{5}$  &  198  &  4.90 $10^{-2}$  &  267  &  4.99 $10^{-2}$\\
2.26 $10^{5}$  &  146  &  4.99 $10^{-2}$  &  272  &  4.95 $10^{-2}$\\
2.06 $10^{5}$  &  188  &  4.99 $10^{-2}$  &  329  &  4.90 $10^{-2}$\\
1.88 $10^{5}$  &  290  &  4.96 $10^{-2}$  &  356  &  4.91 $10^{-2}$\\
1.71 $10^{5}$  &  316  &  4.99 $10^{-2}$  &  317  &  4.85 $10^{-2}$\\
1.55 $10^{5}$  &  326  &  4.98 $10^{-2}$  &  386  &  4.94 $10^{-2}$\\
1.41 $10^{5}$  &  390  &  4.99 $10^{-2}$  &  355  &  4.94 $10^{-2}$\\
1.29 $10^{5}$  &  298  &  5.00 $10^{-2}$  &  376  &  4.97 $10^{-2}$\\
1.17 $10^{5}$  &  329  &  4.88 $10^{-2}$  &  354  &  4.88 $10^{-2}$\\
1.07 $10^{5}$  &  490  &  4.97 $10^{-2}$  &  286  &  4.93 $10^{-2}$\\
9.71 $10^4$  &  615  &  4.91 $10^{-2}$  &  400  &  5.00 $10^{-2}$\\
8.84 $10^4$  &  664  &  4.99 $10^{-2}$  &  504  &  4.95 $10^{-2}$\\
8.05 $10^4$  &  638  &  4.94 $10^{-2}$  &  416  &  4.81 $10^{-2}$\\
7.33 $10^4$  &  436  &  4.94 $10^{-2}$  &  371  &  4.92 $10^{-2}$\\
6.67 $10^4$  &  318  &  4.96 $10^{-2}$  &  319  &  4.70 $10^{-2}$\\
6.07 $10^4$  &  291  &  4.98 $10^{-2}$  &  294  &  4.92 $10^{-2}$\\
5.53 $10^4$  &  351  &  4.85 $10^{-2}$  &  266  &  4.91 $10^{-2}$\\
5.03 $10^4$  &  377  &  4.91 $10^{-2}$  &  258  &  4.78 $10^{-2}$\\
4.58 $10^4$  &  433  &  4.99 $10^{-2}$  &  252  &  4.96 $10^{-2}$\\
4.17 $10^4$  &  480  &  4.97 $10^{-2}$  &  248  &  4.85 $10^{-2}$\\
3.79 $10^4$  &  543  &  4.91 $10^{-2}$  &  247  &  4.75 $10^{-2}$\\
3.45 $10^4$  &  592  &  4.95 $10^{-2}$  &  248  &  5.00 $10^{-2}$\\
3.14 $10^4$  &  547  &  4.99 $10^{-2}$  &  248  &  4.91 $10^{-2}$\\
2.86 $10^4$  &  549  &  4.86 $10^{-2}$  &  248  &  5.00 $10^{-2}$\\
2.61 $10^4$  &  571  &  4.93 $10^{-2}$  &  247  &  4.97 $10^{-2}$\\
2.37 $10^4$  &  618  &  4.99 $10^{-2}$  &  243  &  4.99 $10^{-2}$\\
2.16 $10^4$  &  653  &  4.90 $10^{-2}$  &  244  &  4.99 $10^{-2}$\\
1.97 $10^4$  &  686  &  4.62 $10^{-2}$  &  247  &  4.92 $10^{-2}$\\
1.79 $10^4$  &  723  &  4.83 $10^{-2}$  &  251  &  4.92 $10^{-2}$\\
1.63 $10^4$  &  779  &  4.90 $10^{-2}$  &  254  &  4.97 $10^{-2}$\\
1.48 $10^4$  &  839  &  4.94 $10^{-2}$  &  257  &  4.99 $10^{-2}$\\
1.35 $10^4$  &  855  &  4.90 $10^{-2}$  &  265  &  4.86 $10^{-2}$\\
1.23 $10^4$  &  842  &  4.99 $10^{-2}$  &  267  &  4.95 $10^{-2}$\\
1.12 $10^4$  &  832  &  4.90 $10^{-2}$  &  273  &  4.95 $10^{-2}$\\
1.02 $10^4$  &  832  &  4.94 $10^{-2}$  &  279  &  4.96 $10^{-2}$\\
9.27 $10^3$  &  835  &  4.92 $10^{-2}$  &  280  &  4.98 $10^{-2}$\\
8.44 $10^3$  &  691  &  4.92 $10^{-2}$  &  289  &  4.89 $10^{-2}$\\
7.68 $10^3$  &  843  &  4.96 $10^{-2}$  &  292  &  4.94 $10^{-2}$\\
6.99 $10^3$  &  856  &  4.98 $10^{-2}$  &  296  &  4.97 $10^{-2}$\\
6.36 $10^3$  &  868  &  4.90 $10^{-2}$  &  302  &  5.00 $10^{-2}$\\
5.79 $10^3$  &  875  &  4.93 $10^{-2}$  &  305  &  4.94 $10^{-2}$\\
5.27 $10^3$  &  888  &  4.83 $10^{-2}$  &  311  &  4.95 $10^{-2}$\\
4.80 $10^3$  &  893  &  4.98 $10^{-2}$  &  316  &  4.99 $10^{-2}$\\
4.37 $10^3$  &  907  &  4.94 $10^{-2}$  &  319  &  4.94 $10^{-2}$\\
3.98 $10^3$  &  913  &  4.92 $10^{-2}$  &  327  &  4.94 $10^{-2}$\\
3.62 $10^3$  &  923  &  4.97 $10^{-2}$  &  328  &  4.99 $10^{-2}$\\
3.30 $10^3$  &  944  &  4.99 $10^{-2}$  &  336  &  4.91 $10^{-2}$\\
3.00 $10^3$  &  951  &  4.85 $10^{-2}$  &  339  &  4.94 $10^{-2}$\label{res}\\ 
\end{longtable}\vspace{-0.5cm}
\begin{center}{Table \ref{res}: Iteration table.}\end{center}

\section{Conclusion}
The goal of the study was to allow the computing of micro-wave susceptibility of ferromagnetic particles with thin details.
This last point called for very large meshes (about 300000 degrees of freedom) for which the classical inversion methods
with no preconditionning did not work at all, or required such a large amount of iterations that the computation times for an acceptable 
range of frequencies was far from useful. The strategy presented in this article is an industrial 
computations approach, and obtains interesting results for a large spectrum of benchmark. Computations of realistic experiments have been 
performed (see \cite{Vaast:Numerical,Boust:High,Boust:Micromagnetic}) for physical systems where it was possible to compare
results with physical experiments. Some problems remain, in particular, the strategy developed aims at the laplacian part of the 
total magnetic field whereas some systems are revealed to be principally influenced by the demagnetising field. The next step is
to extend the strategy of the paper in order to include the demagnetisation part of the magnetic field in the approximated
preconditioner. The main difficulty of the extension is algorithmic: to build a good circulant approximation of block Toeplitz
matrices. An another interresting point to study would be the implementation of an efficient parallelised FFT algorithm for distributed
memory systems. The main problem of such an implementation would be the optimisation of the transposition phase of the data
through the memory nodes of the distributed system.

\nocite{*}
\bibliographystyle{unsrt}
\bibliography{aescas1.bib}

\end{document}